\documentclass{scrartcl}

\usepackage[utf8]{inputenc} 
\usepackage[T1]{fontenc}  
\usepackage[english]{babel}
\usepackage{amsmath} 
\usepackage{amsfonts} 
\usepackage{amssymb}
\usepackage{mathrsfs}
\usepackage{authblk}
\usepackage{cite}
\usepackage{typearea}                
\usepackage{scrpage2} 
\usepackage{todonotes}
\usepackage{tikz,pgffor}
\usetikzlibrary{shadows}
\usepackage{mathtools}
\usepackage{tabularx}
\usepackage{rotating,graphics, graphicx, wrapfig}
\usepackage{pgffor}
\usepackage{a4wide}
\addtokomafont{pagehead}{\scshape}

\usetikzlibrary{decorations.pathreplacing}
\usetikzlibrary{topaths,calc}
\usepackage{float}
\usetikzlibrary{snakes,arrows,patterns}

\newtheorem{defsatzusw}{}[section]
\newtheorem{definition}[defsatzusw]{Definition}
\newtheorem{theorem}[defsatzusw]{Theorem}

\newtheorem{corollary}[defsatzusw]{Corollary}
\newtheorem{observation}[defsatzusw]{Observation}

\newtheorem{conjecture}[defsatzusw]{Conjecture}

\DeclareMathOperator{\diam}{diam}

\newenvironment{proof}{\paragraph{$\mathit{Proof.}$}
\vspace{-.4cm}\hspace{-.3cm}}{\vspace{-.2cm}\begin{flushright}$\Box$\end{flushright}}

\newenvironment{alternproof}{\paragraph{$\mathit{Alternative~proof.}$}
	\vspace{-.4cm}\hspace{-.3cm}}{\vspace{-.2cm}\begin{flushright}$\Box$\end{flushright}}

\definecolor{mygray}{HTML}{E6E6E6}
\definecolor{mygreen}{HTML}{298A08}
\definecolor{myred}{HTML}{B40404}

\pgfdeclarepatternformonly{newcrosshatch}{\pgfqpoint{-1pt}{-1pt}}{\pgfqpoint{4pt}{4pt}}{\pgfqpoint{3pt}{3pt}}%
{
	\pgfsetstrokeopacity{0.7}
	\pgfsetlinewidth{0.3pt}
	\pgfpathmoveto{\pgfqpoint{3.2pt}{0pt}}
	\pgfpathlineto{\pgfqpoint{0pt}{3.2pt}}
	\pgfpathmoveto{\pgfqpoint{0pt}{0pt}}
	\pgfpathlineto{\pgfqpoint{3.2pt}{3.2pt}}
	\pgfusepath{stroke}
}

\title{\scshape{On the Burning Number of $p$-Caterpillars}}
\author[1]{Michaela Hiller}
\author[1]{Eberhard Triesch}
\author[1]{Arie M.C.A. Koster}
\affil[1]{\small Lehrstuhl II f\"ur Mathematik, RWTH Aachen, 52062 Aachen, Germany }
\date{\vspace{-5ex}}

\begin{document}
\setkomafont{sectioning}{\normalcolor\bfseries}	
\maketitle

\begin{abstract}\noindent
	\textbf{Abstract.}
	The burning number is a recently introduced graph parameter indicating the spreading speed of content in a graph through its edges. 
	While the conjectured upper bound on the necessary numbers of time steps until all vertices are reached is proven for some specific graph classes it remains open for trees in general. We present two different proofs for ordinary caterpillars 
	and prove the conjecture for a generalised version of caterpillars and for trees with a sufficient amount of leaves.
	
	\noindent
	Furthermore, determining the burning number for spider graphs, trees with maximum degree three and path-forests is known to be $\mathcal{NP}$-complete, however, we show that the complexity is already inherent in caterpillars with maximum degree three.
\end{abstract}


\section{Introduction}


Given an undirected graph $G=(V,E)$, the burning number $b(G)$ indicates the minimum number of steps to inflame the whole graph while in each time step the fire spreads from all burning vertices to their neighbours and one additional vertex can be lit.
This concept was introduced as a possible representation of the spread of content in an online social network in \cite{HowToBurnAGraph}, but also other issues, e.g. the contagion of illnesses, can be modelled.

A sequence of vertices $B=(b_1,\dots,b_{m})$ is said to be a burning sequence or burning strategy if the vertices lit successively burn off the whole graph in $m$ steps. For $m=b(G)$, we say $B$ is an optimum burning sequence resp. strategy. The set of all vertices which receive the fire from a vertex $b_i$ (or theoretically would, if they were not already burning) together with $b_i$ itself is called a burning circle and is denoted by $V_i$.
Thus, finding a burning strategy can be reformulated to a covering problem $V= V_1 \cup \dots \cup V_m$.
Obviously, the extend of a burning circle is given by $\diam(V_i)+1=2i-1$.
We denote the problem of determining the burning number for a graph by \textsc{Burning Number}.

In 2014, an upper bound for the burning number was conjectured for all connected graphs \cite{HowToBurnAGraph}.
\begin{conjecture}[Burning Number Conjecture]
	For any connected graph $G$ of order $n$ it holds $b(G)\leq \left\lceil\sqrt{n}\right\rceil$.
\end{conjecture}
\noindent
The conjecture is proven for paths, cycles, Hamiltonian graphs and spiders \cite{SpidersAndPathForests}. Further, it can easily be checked that graphs with a small vertex number fulfil the conjecture.
For paths whose length is a square number the conjecture holds with equality and, as shown in \cite{HowToBurnAGraph}, the conjecture is true for all connected graphs if it holds for trees in general.

Firstly, in Section \ref{sec:cater} the Burning Number Conjecture is proven for caterpillars in two different ways: once by using the principle of infinite descent and alternatively algorithmically, yielding a burning strategy complying with the conjectured bound.
Subsequently,  in Section \ref{sec:NP_cater} we show that \textsc{Burning Number} is $\mathcal{NP}$-complete for caterpillars. In Section \ref{sec:p-cater} we focus on the validity of the conjecture for $2$-caterpillars and $p$-caterpillars with a sufficient amount of leaves relative to the order of the graph.\\


\section{The Burning Number Conjecture for Caterpillars}\label{sec:cater}

\begin{figure}[H]
	\begin{center}
		\begin{tikzpicture}[decoration=brace,node distance=5em, every node/.style={scale=0.6},scale=.7]
		\node(1) at (-0.866-5,-.5)[circle,fill= black]{};
		\node(2) at (0.866-5,-.5)[circle,fill= black] {};
		\node(3) at (0-5,1)[circle,fill= black] {};
		\draw[line width=0.75pt, black](1)to(2);
		\draw[line width=0.75pt, black](3)to(2);
		\draw[line width=0.75pt, black](1)to(3);

		\node(4) at (0,0)[circle,fill= black] {};
		\node(5) at (0,1)[circle,fill= black] {};
		\node(6) at (0,2)[circle,fill= black] {};
		\node(7) at (-0.866,-.5)[circle,fill= black] {};
		\node(8) at (2*-0.866,2*-.5)[circle,fill= black] {};
		\node(9) at (0.866,-.5)[circle,fill= black] {};
		\node(10) at (2*0.866,2*-.5)[circle,fill= black] {};
		\draw[line width=0.75pt, black](4)to(10);
		\draw[line width=0.75pt, black](4)to(8);
		\draw[line width=0.75pt, black](4)to(6);
		\end{tikzpicture}
		\caption{Forbidden minors $C_3$ and $S_{2,2,2}$ in a caterpillar.}
		\label{minors}		
	\end{center}
\end{figure}
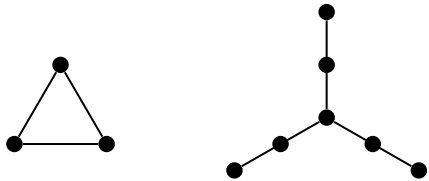

In this section we investigate the Burning Number Conjecture for caterpillars, trees in which all vertices are within the distance one of a central spine
or more vivid:
\begin{quote} 
	`A caterpillar is a tree which metamorphoses into a path when its cocoon of endpoints is removed.' \cite{TheNumberOfCaterpillars}
\end{quote}
\noindent
Consequently, the graph class of caterpillars can also be described by forbidden minors $C_3$ and $S_{2,2,2}$ as in Figure \ref{minors}.

Let $G=(V,E)$ denote a caterpillar with $n\coloneqq|V|$ vertices, a spine $P_l=\{v_1,\dots,v_l\}$ of length $l$ and $n-l$ vertices adjacent to $P_l\setminus\{v_1,v_l\}$, which we call legs. 
We assume $l\geq 4$ and $n\geq l+2$; otherwise $G$ is a spider graph and the conjecture holds. Further, it can easily be seen that the conjecture is true for all graphs with $n\leq 9$.

\noindent
Applying the (proven) conjecture for paths to the Spine $P_l$, we clearly get the following upper bound for the caterpillar.

\begin{observation}
	For a caterpillar $G$ it holds $b(G)\leq \left\lceil\sqrt{l}\right\rceil +1$. 
	Thus, for $\left\lceil\sqrt{n}\right\rceil \geq \Big\lceil\sqrt{l}\Big\rceil +1$ the conjecture is proven to be true.
\end{observation}

\noindent
In fact, the conjecture holds for all caterpillars, which can be shown using the principle of infinite descent and some number-theoretical considerations.

\begin{theorem}[Burning Number Conjecture for Caterpillars]\label{Conjecture1Caterpillar}
	The burning number of a caterpillar $G$ satisfies $b(G)\leq \left\lceil\sqrt{n}\right\rceil$.
	
	\begin{proof}
		Let the graph $G$ be a caterpillar and a minimum counterexample regarding $n$ with $b(G) > \left\lceil\sqrt{n}\right\rceil\eqqcolon k$. We distinguish two cases:
		\begin{enumerate}
			\item[] If either the spine vertex $v_{2k-1}$ does not have any legs or $v_{2k-1}$ has a leg, but at least one of the vertices $v_1,\dots,v_{2k-2}$ has an adjacent leg as well, we remove the largest burning circle $V_1$ with extend $\diam(V_1)+1=2k-1$ without loss of generality at the end of the spine $P_l$. Depending on whether $v_{2k-1}$ is legless or not we shorten the spine by $(2k-1)$ resp. $(2k-1)-1$ vertices to maintain the connectivity.\\
			
			\begin{minipage}{\linewidth}
				\centering
					\begin{tikzpicture}[decoration=brace,node distance=5em, every node/.style={scale=0.6},scale=0.7]
					\node(1) at (0,0)[circle,fill=black] {};
					\node(2) at (-1,0)[circle,fill=black] {};
					\node(3) at (-2,0)[circle,fill=black] {};
					\node(4) at (-3,0)[circle,fill=gray] {};
					\node(5) at (-4,0)[circle,fill= gray] {};
					\node(6) at (-5,0) [circle,fill= gray] {};
					\node(7) at (-2,1)[circle,fill=black] {};
					\node(8) at (-2,-1)[circle,fill=black] {};
					\node(9) at (-3,1)[circle,fill=gray] {};
					\node(10) at (-3,-1)[circle,fill=gray] {};
					
					\node(11) at (1,0) {};
					\node(12) at (0,1){};
					\node(13) at (0,-1) {};
					\node(14) at (-1,1) {};
					\node(15) at (-1,-1) {};
					\node(16) at (-4,1) {};
					\node(17) at (-4,-1) {};
					\node(18) at (-5,1){};
					\node(19) at (-5,-1){};
					\node(20) at (-6,0) {};
					
					\draw[line width=.75pt](1)to(2);
					\draw[line width=.75pt](2)to(3);
					\draw[line width=.75pt, gray](3)to(4);
					\draw[line width=.75pt, gray](4)to(5);
					\draw[line width=.75pt, gray](5)to(6);
					\draw[line width=.75pt](3)to(7);
					\draw[line width=.75pt](3)to(8);
					\draw[line width=.75pt, gray](4)to(9);
					\draw[line width=.75pt, gray](4)to(10);
					\draw[line width=1.5pt, loosely dotted](11)to(1);
					\draw[line width=1.5pt, loosely dotted](12)to(1);
					\draw[line width=1.5pt, loosely dotted](13)to(1);
					\draw[line width=1.5pt, loosely dotted](14)to(2);
					\draw[line width=1.5pt, loosely dotted](15)to(2);
					\draw[line width=1.5pt, loosely dotted, gray](16)to(5);
					\draw[line width=1.5pt, loosely dotted, gray](17)to(5);
					\draw[line width=1.5pt, loosely dotted, gray](18)to(6);
					\draw[line width=1.5pt, loosely dotted, gray](19)to(6);
					\draw[line width=1.5pt, loosely dotted, gray](20)to(6);
					
					\draw[draw=gray, fill=gray, opacity=0.2] (-4.15,0) ellipse (2.4cm and 1.4cm);
					
					\draw[decorate, yshift=-4ex] (-1.75,-.8) -- node[below=0.5ex] {~ } (-6.4,-.8);
					\node at (-4.075,-2) [scale=1.5]{$2k-1$};
					\node at (-6,1.3) [scale=1.5]{$V_1$};
					\end{tikzpicture}
					\label{Bild1}
					
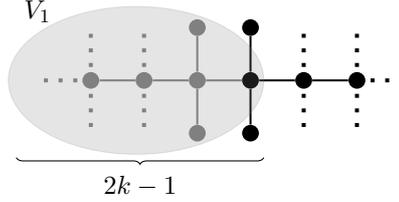
\captionof{figure}{We generate $G^\prime$ by removing the grey vertices of $V_1$ in the minimum counterexample $G$.}
				\end{minipage}
				
			\noindent
			In both sub-cases, we obtain a new caterpillar $G^\prime$ with 
			\begin{align*}
				l^\prime &\leq l-(2k-1)+2 &&\hspace{-3.5cm}= l-2\left\lceil\sqrt{n}\right\rceil+3,\\
				n^\prime &\leq n-(2k-1) &&\hspace{-3.5cm}= n-2\left\lceil\sqrt{n}\right\rceil+1, 
			\end{align*}
			and for the burning number of $G^\prime$ it follows $b(G^\prime) > \left\lceil\sqrt{n}\right\rceil-1$, otherwise $G$ would not be a counterexample. 
			Since $G$ is minimum, it further holds $b(G^\prime) \leq \left\lceil\sqrt{n^\prime}\right\rceil$. This yields
			\begin{align*}
				&&\hspace{-3cm}\left\lceil\sqrt{n^\prime}\right\rceil \geq b(G^\prime) & \;>\; \left\lceil\sqrt{n}\right\rceil-1 
			\end{align*}
			and thus, $\left\lceil\sqrt{n^\prime}\right\rceil= \left\lceil\sqrt{n}\right\rceil$. With the estimate from above we get $\left\lceil\sqrt{n-2\left\lceil\sqrt{n}\right\rceil+1}\right\rceil =\left\lceil\sqrt{n}\right\rceil$ and therefore, the two radicands lie between the same square numbers $\left\lceil\sqrt{n}\right\rceil^2$ and $\left(\left\lceil\sqrt{n}\right\rceil-1\right)^2$. As a consequence
			\begin{align*}
				n-\left(\left\lceil\sqrt{n}\right\rceil-1\right)^2 &\; \geq\; 2\left\lceil\sqrt{n}\right\rceil-1+1,
			\end{align*}
			or equivalently, $n  \geq \left(\left\lceil\sqrt{n}\right\rceil-1\right)^2+2\left\lceil\sqrt{n}\right\rceil= \left\lceil\sqrt{n}\right\rceil^2+1$. This is a contradiction.
			
			\item[] If otherwise $v_1,\dots,v_{2k-2}$ are legless but $v_{2k-1}$ is not, we remove the two largest burning circles $V_1$ with extend $\diam(V_1)+1=2k-1$ and $V_2$ with extend $\diam(V_2)+1=2k-3$ without loss of generality at the end of the spine $P_l$. We shorten the spine by $(2k-3)+(2k-1)-1$ vertices.\\
					
			\begin{minipage}{\linewidth}
				\centering
				\begin{tikzpicture}[decoration=brace,node distance=5em, every node/.style={scale=0.6},scale=0.7]
				\node(1) at (-1.5,0)[circle,fill=gray]{};
				\node(2) at (1,0)[circle,fill=gray] {};
				\node(3) at (2,0)[circle,fill=gray] {};
				\node(4) at (3.5,0)[circle,fill=gray] {};
				\node(5) at (4.5,0)[circle,fill=black] {};
				\node(6) at (5.5,0) [circle,fill= black] {};
				\node(7) at (2,1)[circle,fill=gray]{};
				\node(8) at (2,-1) {};
				\node(9) at (3.5,1){};
				\node(10) at (3.5,-1){};
				
				\node(11) at (-2.5,0) [circle,fill=gray] {};
				\node(12) at (0,0)[circle,fill=gray]{};
				\node(16) at (4.5,1) {};
				\node(17) at (4.5,-1) {};
				\node(18) at (5.5,1){};
				\node(19) at (5.5,-1){};
				\node(20) at (6.5,0)[circle,fill= black]{};
				\node(21) at (6.5,1){};
				\node(22) at (6.5,-1){};
				\node(23) at (7.5,0) {};
				
				\draw[line width=1.5pt, loosely dotted, gray](1)to(12);
				\draw[line width=.75pt, gray](2)to(3);
				\draw[line width=1.5pt, loosely dotted, gray](3)to(4);
				\draw[line width=.75pt,gray](4)to(5);
				\draw[line width=.75pt](5)to(6);
				\draw[line width=.75pt, gray](3)to(7);
				\draw[line width=1.5pt, loosely dotted, gray](3)to(8);
				\draw[line width=1.5pt, loosely dotted,gray](4)to(9);
				\draw[line width=1.5pt, loosely dotted,gray](4)to(10);
				\draw[line width=.75pt, gray](11)to(1);
				\draw[line width=.75pt, gray](12)to(2);
				\draw[line width=1.5pt, loosely dotted,black](16)to(5);
				\draw[line width=1.5pt, loosely dotted,black](17)to(5);
				\draw[line width=1.5pt, loosely dotted](18)to(6);
				\draw[line width=1.5pt, loosely dotted](19)to(6);
				\draw[line width=.75pt](20)to(6);
				\draw[line width=1.5pt, loosely dotted](20)to(21);
				\draw[line width=1.5pt, loosely dotted](20)to(22);
				\draw[line width=1.5pt, loosely dotted](20)to(23);
				
				\draw[draw=gray, fill=gray, opacity=0.2] (-1.15,0) ellipse (1.65cm and 1.2cm);
				\draw[draw=gray, fill=gray, opacity=0.2] (2.84,0) ellipse (2.3cm and 1.3cm);
				
				
				\node at (-2.5,-.4) [scale=1.25]{$v_1$};
				\node at (-1.5,-.4) [scale=1.25]{$v_2$};
				\node at (-.1,-.4) [scale=1.25]{$v_{2k-3}$};
				\node at (1.2,.35) [scale=1.25]{$v_{2k-2}$};
				\node at (1.9,-.4) [scale=1.25, fill=mygray]{\phantom{..}};
				\node at (1.9,-.4) [scale=1.25]{$v_{2k-1}$};
				\node at (3.35,-.4) [scale=1.25, fill=mygray]{\phantom{..}};
				\node at (3.35,-.4) [scale=1.25]{$v_{4k-5}$};
				\node at (4.5,-.4) [scale=1.25, fill=mygray]{\phantom{..}};
				\node at (4.5,-.4) [scale=1.25]{$v_{4k-4}$};
				
				\draw[decorate, yshift=-4ex] (.4,-.8) -- node[below=0.5ex] {~} (-2.65,-.8);
				\draw[decorate, yshift=-4ex] (5.,-.8) -- node[below=0.5ex] {~} (.6,-.8);
				\node at (-1.125,-2) [scale=1.5]{$2k-3$};
				\node at (2.8,-2) [scale=1.5]{$2k-1$};
				\node at (1.25,1.3) [scale=1.5]{$V_1$};
				\node at (-2.3,1.3) [scale=1.5]{$V_2$};
				\end{tikzpicture}
				\label{Bild4}
				
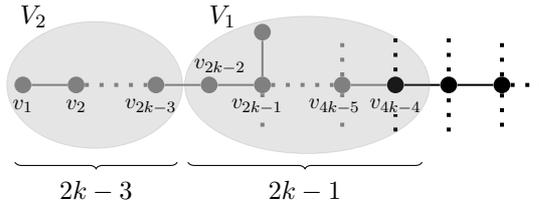
\captionof{figure}{Let $v_{2k-1}$ have adjacent legs and $v_1,\dots,v_{2k-2}$ be legless. We remove the grey vertices.}
			\end{minipage}
		
			\noindent
			Analogously to the first case, for the remaining caterpillar $G^{\prime\prime}$ it follows
			\begin{align*}
				l^{\prime\prime} &\leq l-(2k-3)-(2k-1)+2 ,\\
				n^{\prime\prime} &\leq n-(2k-3)-(2k-1)+1-1 \leq \left( \left\lceil\sqrt{n}\right\rceil-2\right)^2 , 
			\end{align*}
			and $b(G^{\prime\prime}) > \left\lceil\sqrt{n}\right\rceil-2$, otherwise $G$ would not be a counterexample. 
			Since $G$ is minimum, it further holds $b(G^\prime) \leq \left\lceil\sqrt{n^{\prime\prime}}\right\rceil$. This yields the contradiction
			$$\left\lceil\sqrt{n}\right\rceil-2 \;< \;b(G^{\prime\prime})\;\leq\; \left\lceil\sqrt{n^{\prime\prime}}\right\rceil \;\leq \;\left\lceil\sqrt{n}\right\rceil-2.$$
		\end{enumerate}

\noindent
		Therefore, the minimum counterexample cannot exist and the theorem holds true.\vspace{-.5cm}
	\end{proof}
\vspace{.5cm}
\noindent
\textnormal{The following alternative proof works without the principle of infinite descent and provides a burning strategy in $\left\lceil\sqrt{n}\right\rceil$ steps for all caterpillars.}

	\begin{alternproof}
		Let again $k\coloneqq\left\lceil\sqrt{n}\right\rceil$ denote the 
		maximum number of steps such that the conjecture still holds.
		Recursively removing burning circles to reduce the vertex number at least down to the next smaller square number, we consider two cases: \begin{enumerate}
			\item[] In the first case, $v_{2k-1}\in P_l$ does not have any legs. After deleting $v_1,\dots,v_{2k-1}$ with all adjacent legs the remaining graph has at most $n-(2k-1)\leq \left\lceil\sqrt{n}\right\rceil^2-2\left\lceil\sqrt{n}\right\rceil+1=(\left\lceil\sqrt{n}\right\rceil-1)^2$ vertices.
			
			\begin{minipage}{\linewidth}
				\centering
					\begin{tikzpicture}[decoration=brace,node distance=5em, every node/.style={scale=0.6},scale=0.7]
					\node(1) at (-.5,0)[circle,fill=gray]{};
					\node(2) at (1,0)[circle,fill=gray] {};
					\node(3) at (2,0)[circle,fill=gray] {};
					\node(4) at (3,0)[circle,fill=black] {};
					\node(5) at (4,0)[circle,fill=black] {};
					\node(6) at (5,0) [circle,fill= black] {};
					\node(7) at (2,1){};
					\node(8) at (2,-1) {};
					\node(9) at (3,1){};
					\node(10) at (3,-1){};
					
					\node(11) at (-1.5,0) [circle,fill=gray] {};
					\node(12) at (-.5,1){};
					\node(13) at (-.5,-1) {};
					\node(14) at (1,1) {};
					\node(15) at (1,-1) {};
					\node(16) at (4,1) {};
					\node(17) at (4,-1) {};
					\node(18) at (5,1){};
					\node(19) at (5,-1){};
					\node(20) at (6,0) {};
					
					\draw[line width=1.5pt, loosely dotted, gray](1)to(2);
					\draw[line width=.75pt, gray](2)to(3);
					\draw[line width=.75pt,gray](3)to(4);
					\draw[line width=.75pt](4)to(5);
					\draw[line width=.75pt](5)to(6);
					\draw[line width=1.5pt, loosely dotted](4)to(9);
					\draw[line width=1.5pt, loosely dotted](4)to(10);
					\draw[line width=.75pt, gray](11)to(1);
					\draw[line width=1.5pt, loosely dotted, gray](12)to(1);
					\draw[line width=1.5pt, loosely dotted, gray](13)to(1);
					\draw[line width=1.5pt, loosely dotted, gray](14)to(2);
					\draw[line width=1.5pt, loosely dotted, gray](15)to(2);
					\draw[line width=1.5pt, loosely dotted](16)to(5);
					\draw[line width=1.5pt, loosely dotted](17)to(5);
					\draw[line width=1.5pt, loosely dotted](18)to(6);
					\draw[line width=1.5pt, loosely dotted](19)to(6);
					\draw[line width=1.5pt, loosely dotted](20)to(6);
					
					\draw[draw=gray, fill=gray, opacity=0.2] (0.3,0) ellipse (2.3cm and 1.4cm);
					
					
					\node at (-1.5,-.415) [scale=1.25]{$v_1$};
					\node at (-.5,-.415) [scale=1.25, fill=mygray]{$v_2$};
					\node at (.8,-.415) [scale=1.25, fill=mygray]{$v_{2k-2}$};
					\node at (2,-.415) [scale=1.25]{$v_{2k-1}$};
					
					\draw[decorate, yshift=-4ex] (2.4,-.8) -- node[below=0.5ex] {~ } (-1.7,-.8);
					\node at (0.35,-2) [scale=1.5]{$2k-1$};
					\node at (2,1.3) [scale=1.5]{$V_1$};
					\end{tikzpicture}
					\label{Bild2}
				
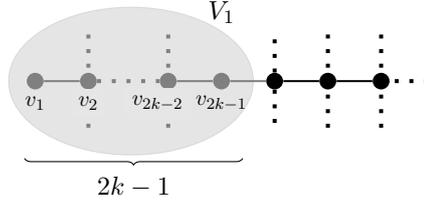
\captionof{figure}{In the first case, $v_{2k-1}$ is legless and we delete the grey vertices.}
			\end{minipage}
	
			\item[]In the other case, we distinguish whether any of the vertices $v_1,\dots,v_{2k-2}$ has an adjacent leg or not. 
			If not all of these spine vertices are legless, we remove $v_1,\dots,v_{2k-2}$ together with their legs. 
			Again the vertex set of the remaining graph contains - just as in the first case - at most $(\left\lceil\sqrt{n}\right\rceil-1)^2$ vertices.\\
			
			 \begin{minipage}{\linewidth}
				\centering
					\begin{tikzpicture}[decoration=brace,node distance=5em, every node/.style={scale=0.6},scale=0.7]
					\node(1) at (-.5,0)[circle,fill=gray]{};
					\node(2) at (1,0)[circle,fill=gray] {};
					\node(3) at (2,0)[circle,fill=black] {};
					\node(4) at (3,0)[circle,fill=black] {};
					\node(5) at (4,0)[circle,fill=black] {};
					\node(6) at (5,0) [circle,fill= black] {};
					\node(7) at (2,1)[circle,fill=black]{};
					\node(8) at (2,-1) {};
					\node(9) at (3,1){};
					\node(10) at (3,-1){};
					
					\node(11) at (-1.5,0) [circle,fill=gray] {};
					\node(12) at (-.5,1){};
					\node(13) at (-.5,-1) {};
					\node(14) at (1,1) {};
					\node(15) at (1,-1) {};
					\node(16) at (4,1) {};
					\node(17) at (4,-1) {};
					\node(18) at (5,1){};
					\node(19) at (5,-1){};
					\node(20) at (6,0) {};
					
					\draw[line width=1.5pt, loosely dotted, gray](1)to(2);
					\draw[line width=.75pt, gray](2)to(3);
					\draw[line width=.75pt](3)to(4);
					\draw[line width=.75pt](4)to(5);
					\draw[line width=.75pt](5)to(6);
					\draw[line width=.75pt](3)to(7);
					\draw[line width=1.5pt, loosely dotted](3)to(8);
					\draw[line width=1.5pt, loosely dotted](4)to(9);
					\draw[line width=1.5pt, loosely dotted](4)to(10);
					\draw[line width=.75pt, gray](11)to(1);
					\draw[line width=1.5pt, loosely dotted, gray](12)to(1);
					\draw[line width=1.5pt, loosely dotted, gray](13)to(1);
					\draw[line width=1.5pt, loosely dotted, gray](14)to(2);
					\draw[line width=1.5pt, loosely dotted, gray](15)to(2);
					\draw[line width=1.5pt, loosely dotted](16)to(5);
					\draw[line width=1.5pt, loosely dotted](17)to(5);
					\draw[line width=1.5pt, loosely dotted](18)to(6);
					\draw[line width=1.5pt, loosely dotted](19)to(6);
					\draw[line width=1.5pt, loosely dotted](20)to(6);
					
					\draw[draw=gray, fill=gray, opacity=0.2] (0.3,-.1) ellipse (2.3cm and 1.2cm);
					
					
					\node at (-1.5,-.4) [scale=1.25]{$v_1$};
					\node at (-.5,-.4) [scale=1.25,fill=mygray]{\phantom{..}};
					\node at (-.5,-.4) [scale=1.25]{$v_2$};
					\node at (.8,-.4) [scale=1.25,fill=mygray]{\phantom{..}};
					\node at (.8,-.4) [scale=1.25]{$v_{2k-2}$};
					\node at (2.,-.4) [scale=1.25, fill=mygray]{\phantom{-}};
					\node at (2.,-.4) [scale=1.25]{$v_{2k-1}$};
					
					\draw[decorate, yshift=-4ex] (2.4,-.8) -- node[below=0.5ex] {~ } (-1.7,-.8);
					\node at (.35,-2) [scale=1.5]{$2k-1$};
					\node at (-1.5,1.1) [scale=1.5]{$V_1$};
					\end{tikzpicture}

					\label{Bild3}
				
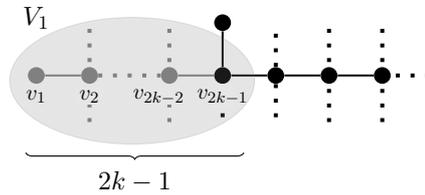
\captionof{figure}{We assume $v_{2k-1}$ and at least one of $v_1,\dots,v_{2k-2}$ to have legs and delete the grey vertices.}
			\end{minipage}
		
			Otherwise, if $v_1,\dots,v_{2k-2}$ are legless and $v_{2k-1}$ has an adjacent leg, we delete $v_1,\dots,v_{2k-3}$ and further $v_{(2k-3)+1},\dots,v_{(2k-3)+(2k-2)}$ with all their legs (at least the leg adjazent to $v_{2k-1}$) such that the new graph consists of at most $n-(2k-3)-(2k-2)-1\leq \left\lceil\sqrt{n}\right\rceil^2-(2\left\lceil\sqrt{n}\right\rceil-1)-(2\left\lceil\sqrt{n}\right\rceil-3)=(\left\lceil\sqrt{n}\right\rceil-2)^2$ vertices.
		\end{enumerate}
	Hence, after the vertex removal the order of the remaining graph $G^\prime$ decreases at least to $n^\prime\leq (\left\lceil\sqrt{n}\right\rceil-1)^2$ and the claim follows recursively.
	\vspace{-.53cm}
	\end{alternproof}
\end{theorem}\vspace{.25cm}
\noindent
It can easily be seen that the alternative proof yields an algorithm to burn a caterpillar in $\left\lceil\sqrt{n}\right\rceil$ steps, though may not necessarily be optimum.\\

\section{The $\mathcal{NP}$-Completeness of the Burning Number Problem for Caterpillars}\label{sec:NP_cater}


The $\mathcal{NP}$-completeness of determining the burning number for caterpillars indicates the unstructured nature of the problem, as the difficulty or complexity is already hidden in such a simple graph class. 
Our proof is structured similar to the proof for trees of maximum degree three in \cite{BurningAGraphIsHard} and uses a reduction from \textsc{Distinct $3$-Partition}.\\

\noindent
\textbf{Problem:} \textsc{Distinct $3$-Partition}\\
\textbf{Instance:}  \begin{minipage}[t]{\linewidth-2cm}
						A set $X=\{a_1, \dots, a_{3n}\}$ of $3n$ distinct positive integers and a positive integer $S$, fulfilling $\sum_{i=1}^{3n}a_i=n\cdot S$ with $\frac{S}{4}<a_i<\frac{S}{2}$ for all $1\leq i\leq 3n$.
					\end{minipage}
\textbf{Question:} 	\begin{minipage}[t]{\linewidth-2cm}
						Can $X$ be partitioned into $n$ triples each of whose elements sum up to $S$? \\
					\end{minipage}
\textsc{Distinct $3$-Partition} is $\mathcal{NP}$-complete in the strong sense as shown in \cite{MultigraphRealizationsOfDegreeSequences}, which means the problem remains $\mathcal{NP}$-complete, even if $S$ is bounded from above by a polynomial in $n$.

\begin{theorem}
	\textsc{Burning Number} is $\mathcal{NP}$-complete for caterpillars of maximum degree three.
	
	\begin{proof}
		\textsc{Burning Number} is in $\mathcal{NP}$, as a burning sequence for a graph can be verified in polynomial time by checking whether the whole vertex set is covered by the union of the corresponding burning circles.
		
		To prove the $\mathcal{NP}$-completeness, we reduce \textsc{Distinct $3$-Partition} in polynomial time to \textsc{Burning Number}. 
		Given an instance for \textsc{Distinct $3$-Partition} as stated above, we denote $m\coloneqq\max\{a_i~|~a_i\in X\}$, $\underline{m}\coloneqq\{1,\dots,m\}$ and $Y\coloneqq\underline{m}\setminus X$. 
		Transferred to the universe of \textsc{Burning Number}, we get $X'\coloneqq\{2a_i-1~|~a_i\in X\}$, $S'\coloneqq2S-3$, $\mathcal{O}_m\coloneqq\{2i-1~|~i\in \underline{m}\}$ and $Y'\coloneqq\mathcal{O}_m\setminus X'$.
		
		Now, we construct a caterpillar $G$ of maximum degree three as follows: For each triple whose unknown elements should add up to $S$ we build a path $Q^{X'}_i$ (for all $1\leq i \leq n$) of order $S'$ and for all numbers in $Y$ (which are not available for the triples) a separate path $Q^{Y'}_i$ (for all $1\leq i \leq m-3n$) of order $Y'$. 
		The resulting path forest 
		$$\bigcup\limits_{i=1}^n Q_i^{X'}\cup \bigcup\limits_{i=1}^{m-3n} Q_i^{Y'}$$
		corresponds to $\bigcup_{i=1}^m P_{2i-1}$ and thus, it can be burnt in $m$ steps. 
		Next, we need to connect the graph by using caterpillars to keep the individual paths isolated. 
		In order to do so, we need at most $m+1$ caterpillars $G_1, \dots, G_{m+1}$, whereby $G_i$ has a spine of length $2(2m+1-i)+1$ with exactly one leg attached to each spine vertex (except the two terminal vertices). 
		The caterpillars and the paths are arranged alternately until only caterpillars are left, which are then placed at the end. The subgraphs are connected through an edge between their end vertices.
		We denote the longest path in $G$ by $P_l$ and get 
		\begin{align*}
			l&=\Bigg|\bigcup\limits_{i=1}^n V\Big(Q_i^{X'}\Big)\Bigg|\cup \Bigg|\bigcup\limits_{i=1}^{m-3n} V\Big(Q_i^{Y'}\Big)\Bigg|\cup \Bigg|\bigcup\limits_{i=1}^{m+1} V\big(P_{2(2m+1-i)+1}\big)\Bigg|\\
			   &=\sum\limits_{i=1}^{m}(2i-1)+ \sum\limits_{i=1}^{m+1}(2(2m+1-i)+1) \\
			   &=\sum\limits_{i=1}^{m}(2i-1)+ \sum\limits_{i=m+1}^{2m+1}(2i-1) \\
			   &= (2m+1)^2.
		\end{align*}
		The inequality in the conjecture is sharp for paths and thus, $b(G)\geq b(P_l) = \left\lceil\sqrt{l}\right\rceil=2m+1$.
		Due to the strong $\mathcal{NP}$-completeness of \textsc{Distinct $3$-Partition}, we can assume $S$ to be in $\mathcal{O}\big(n^{\mathcal{O}(1)}\big)$ and as $m$ is bounded by $S$, the caterpillar $G$ is computed in polynomial time with regard to the input length.
		Further, we constructed the caterpillar $G$ in such a way that, if $X$ can be partitioned into $n$ triples, each of whose elements add up to $S$ (and equivalently $Q_1^{X'},\dots,Q_n^{X'}$ can be partitioned in paths $\{P_i~|~i\in X'\}$), lighting the central spine vertex of caterpillar $G_i$ in step $i$ (for $1\leq i\leq m+1$) and lighting the central vertex of path $P_{2(2m+1-i)+1}$ in step $i$ (for $m+2\leq i\leq 2m+1$) burns the whole graph in $2m+1$ steps. Consequently, it also holds $b(G)\leq 2m+1$ and altogether, $b(G)=2m+1$.
		
		To prove the opposite direction, we assume $b(G)=2m+1$ and let $(x_1, \dots, x_{2m+1})$ be an optimal burning sequence for the caterpillar $G$.
		First, we can observe that $x_i$ is a spine vertex for all $1\leq i\leq 2m+1$ and the burning circles have to be pairwise disjoint, as $l$ is a square number and $b(P_l)=\left\lceil\sqrt{l}\right\rceil$.
		Next, the largest burning circles has to cover $G_1$ with spine $P_{2(2m+1)-1}$. 
		Otherwise, at least two burning circles are needed which would have to intersect at two spine vertices to cover all legs as pictured in Figure \ref{BildNP}. 
		Inductively, $G_i$ has to be covered with the $i$-th largest burning circle, and thus the central spine vertex of $G_i$ has to be lit in the $i$-th step for all $1\leq i\leq m+1$.\\
		
		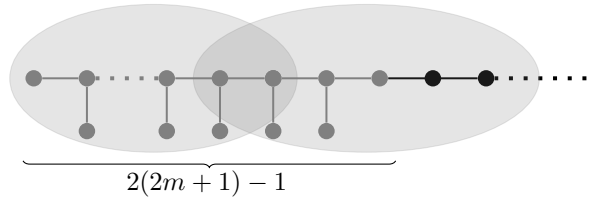
\begin{figure}[H]
		\begin{center}
			\begin{tikzpicture}[decoration=brace,node distance=5em, every node/.style={scale=0.6},scale=0.7]
			\node(1) at (-1.5,0)[circle,fill=gray]{};
			\node(2) at (0,0)[circle,fill=gray] {};
			\node(3) at (2,0)[circle,fill=gray] {};
			\node(4) at (3,0)[circle,fill=gray] {};
			\node(5) at (4,0)[circle,fill=gray] {};
			\node(6) at (5,0) [circle,fill= black] {};
			\node(7) at (1,0)[circle,fill=gray]{};
			\node(8) at (1,-1)[circle,fill=gray]{};
			\node(9) at (2,0)[circle,fill=gray]{};
			\node(10) at (2,-1)[circle,fill=gray]{};
			\node(11) at (-2.5,0) [circle,fill=gray] {};
			\node(12) at (3,-1) [circle,fill=gray] {};
			\node(13) at (-1.5,-1)[circle,fill=gray] {};
			\node(15) at (0,-1)[circle,fill=gray] {};
			\node(16) at (6,0) [circle,fill= black] {};
			\node(20) at (8,0) {};
			
			\draw[line width=1.5pt, loosely dotted, gray](1)to(2);
			\draw[line width=.75pt, gray](2)to(3);
			\draw[line width=.75pt, gray](7)to(8);
			\draw[line width=.75pt, gray](9)to(10);
			\draw[line width=.75pt, gray](4)to(12);
			\draw[line width=.75pt,gray](3)to(4);
			\draw[line width=.75pt,gray](4)to(5);
			\draw[line width=.75pt](5)to(16);
			\draw[line width=.75pt, gray](11)to(1);
			\draw[line width=.75pt,gray](13)to(1);
			\draw[line width=.75pt, gray](15)to(2);
			\draw[line width=1.5pt, loosely dotted](20)to(16);

			\draw[draw=gray, fill=gray, opacity=0.2] (-0.25,0) ellipse (2.7cm and 1.4cm);
			\draw[draw=gray, fill=gray, opacity=0.2] (3.75,0) ellipse (3.25cm and 1.4cm);
			
			\draw[decorate, yshift=-4ex] (4.3,-.9) -- node[below=0.5ex] {~ } (-2.7,-.9);
			\node at (0.75,-2) [scale=1.5]{$2(2m+1)-1$};
		\end{tikzpicture}
		\caption{If we do not cover $G_1$ with the largest burning circles at least two spine vertices are covered twice.}
		\label{BildNP}		
		\end{center}
		\end{figure}
	\vspace{-.4cm}
		\noindent
		Therefore, $\bigcup_{i=1}^{m+1}G_i$ will be burning after $2m+1$ steps induced by $x_1,\dots, x_{m+1}$ and in the last $m$ time steps $x_{m+2} ,\dots,x_{2m+1}$ have to ignite $\bigcup_{i=1}^n Q_i^{X'}\cup \bigcup_{i=1}^{m-3n} Q_i^{Y'}=\bigcup_{i=1}^m P_{2i-1}$, i.e., the remaining subpaths need to be covered by 
		$$\bigcup\limits_{i=m+2}^{2m+1}\mathcal{N}_{2m+1-i}[x_i].$$ 
		As seen before the burning circles have to be disjoint, and thus $\mathcal{N}_{2m+1-i}[x_i]$ has to cover a path of length $2(2m+1-i)-1$ for $m+2\leq i\leq 2m+1$. 
		Hence, each $Q^{Y'}_i$ for $1\leq i \leq m-3n$  is covered by itself and each $Q^{X'}_i$ for $1\leq i \leq n$ is partitioned in paths of lengths $X'$. 
		Since $\frac{S}{4}<a_i<\frac{S}{2}$ by assumption, each partition consists of three elements in $X'$, which add up to $2S-3$.
		By retranslating this $3$-partition of $X'$ to $X$ we obtain the sought-for partition into $n$ triples each of whose elements sum up to $S$.\vspace{-.55cm}
	\end{proof}
\end{theorem}
\noindent
As caterpillars are exactly the trees of pathwidth one the above theorem provides a statement about the complexity of graphs whose spanning trees are caterpillars.

\begin{corollary}
	\textsc{Burning Number} is $\mathcal{NP}$-complete for graphs of pathwidth one.\\
\end{corollary}


\section{The Burning Number Conjecture for $p$-Caterpillars}\label{sec:p-cater}


In this section we turn the study to the more general case of $p$-caterpillars.
\begin{definition}[$p$-Caterpillar]
	A $p$-caterpillar $G$ is a tree in which all vertices are within a distance $p$ of a central spine $P_l=\{v_1,\dots,v_l\}$, which is the longest path in $G$. \\
	Further, $r$-legs of a given $p$-caterpillar are defined as disjoint subtrees of $G- P_l$ with depth $r-1$, for $r\leq p$, whose roots are in distance one of the spine.
	We denote the maximum length of all legs attached to spine vertex $v_i$ by $\ell_{\max}(v_i)$ and the number of all vertices which are connected to the spine via $v_i$ by $\ell_{\Sigma}(v_i)$. 
\end{definition}
\noindent
Thus, the parameter $p$ indicates the maximum length of the legs and for every tree $T$ there is a $p$ such that $T$ can be regarded as a $p$-caterpillar. Obviously, a $1$-caterpillar denotes a `common' caterpillar.

\begin{observation}
	For a $p$-caterpillar $G$ it holds $b(G)\leq \left\lceil\sqrt{l}\right\rceil +p$. 
	Thus, for $\left\lceil\sqrt{n}\right\rceil \geq \left\lceil\sqrt{l}\right\rceil +p$ the conjecture is proven to be true.
\end{observation}
\noindent
Using a similar idea as in the alternative proof of Theorem \ref{Conjecture1Caterpillar}, we can prove the Burning Number Conjecture for $2$-caterpillars.

\begin{theorem}[Burning Number Conjecture for $2$-Caterpillars]\label{Conjecture2Caterpillar}
	The burning number of a $2$-caterpillar $G$ satisfies $b(G)\leq \left\lceil\sqrt{n}\right\rceil$.
	
	\begin{proof}
		As in the alternative proof of Theorem \ref{Conjecture1Caterpillar} we remove recursively the largest burning circles and thereby intend to reduce the number of vertices to fall below the next smaller square number.
		\noindent
		If $\ell_{\max}(v_{2k-2})\leq 1$ and $\ell_{\max}(v_{2k-1})=0$, we delete the vertices $v_1,\dots, v_{2k-1}$ together with all adjacent legs and obtain a graph whose vertex number is at most $\left\lfloor\sqrt{n}\right\rfloor^2$.
		\noindent
		In the case $\ell_{\max}(v_{2k-2})=2$ or $\ell_{\max}(v_{2k-1})\geq 1$ but $\sum_{i=1}^{2k-3}\ell_\Sigma(v_i)\geq 2$, removing the vertices $v_1,\dots, v_{2k-3}$ with their adjacent legs as depicted in Figure \ref{(3)} suffices to undercut $\left\lfloor\sqrt{n}\right\rfloor^2$ vertices in the remaining graph.
		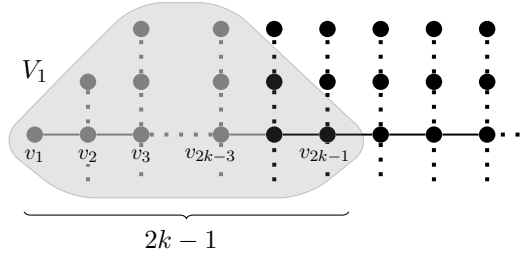
\begin{figure}[H]
			\begin{center}
				\begin{tikzpicture}[decoration=brace,node distance=5em, every node/.style={scale=0.6},scale=0.7]
					\node(1) at (-2.5,0)[circle,fill=gray]{};
					\node(2) at (1,0)[circle,fill=black] {};
					\node(3) at (2,0)[circle,fill=black] {};
					\node(4) at (3,0)[circle,fill=black] {};
					\node(5) at (4,0)[circle,fill=black] {};
					\node(6) at (5,0) [circle,fill= black] {};
					\node(7) at (2,1)[circle,fill=black]{};
					\node(8) at (2,-1) {};
					\node(9) at (3,1)[circle,fill=black]{};
					\node(10) at (3,-1){};
					
					\node(11) at (-3.5,0) [circle,fill=gray] {};
					\node(12) at (-2.5,1)[circle,fill=gray]{};
					\node(13) at (-2.5,-1) {};
					\node(14) at (1,1)[circle,fill=black] {};
					\node(15) at (1,-1) {};
					\node(16) at (4,1) [circle,fill=black]{};
					\node(17) at (4,-1) {};
					\node(18) at (5,1)[circle,fill=black]{};
					\node(19) at (5,-1){};
					\node(20) at (6,0) {};
					\node(21) at (0,0) [circle,fill=gray] {};
					\node(22) at (0,1)[circle,fill=gray]{};
					\node(23) at (0,-1.){};
					\node(24) at (-1.5,0) [circle,fill=gray] {};
					\node(25) at (-1.5,1)[circle,fill=gray]{};
					\node(26) at (-1.5,-1) {};
					
					\node(27) at (-1.5,2)[circle,fill=gray] {};
					\node(28) at (0,2)[circle,fill=gray] {};
					\node(29) at (1,2)[circle,fill=black] {};
					\node(30) at (2,2)[circle,fill=black] {};
					\node(31) at (3,2)[circle,fill=black]{};
					\node(32) at (4,2)[circle,fill=black]{};
					\node(33) at (5,2)[circle,fill=black]{};
					
					\draw[line width=1.5pt, loosely dotted, gray](24)to(21);
					\draw[line width=.75pt](2)to(3);
					\draw[line width=.75pt, gray](2)to(21);
					\draw[line width=.75pt, gray](1)to(24);
					\draw[line width=.75pt](3)to(4);
					\draw[line width=.75pt](4)to(5);
					\draw[line width=.75pt](5)to(6);
					\draw[line width=1.5pt, loosely dotted](3)to(7);
					\draw[line width=1.5pt, loosely dotted](3)to(8);
					\draw[line width=1.5pt, loosely dotted](14)to(29);
					\draw[line width=1.5pt, loosely dotted](7)to(30);
					\draw[line width=1.5pt, loosely dotted](4)to(9);
					\draw[line width=1.5pt, loosely dotted](4)to(10);
					\draw[line width=1.5pt, loosely dotted](9)to(31);
					\draw[line width=1.5pt, loosely dotted](16)to(32);
					\draw[line width=1.5pt, loosely dotted](18)to(33);
					\draw[line width=.75pt, gray](11)to(1);
					\draw[line width=1.5pt, loosely dotted, gray](12)to(1);
					\draw[line width=1.5pt, loosely dotted, gray](13)to(1);
					\draw[line width=1.5pt, loosely dotted](14)to(2);
					\draw[line width=1.5pt, loosely dotted](15)to(2);
					\draw[line width=1.5pt, loosely dotted](16)to(5);
					\draw[line width=1.5pt, loosely dotted](17)to(5);
					\draw[line width=1.5pt, loosely dotted](18)to(6);
					\draw[line width=1.5pt, loosely dotted](19)to(6);
					\draw[line width=1.5pt, loosely dotted](20)to(6);
					\draw[line width=1.5pt, loosely dotted, gray](21)to(22);
					\draw[line width=1.5pt, loosely dotted, gray](21)to(23);
					\draw[line width=1.5pt, loosely dotted, gray](24)to(25);
					\draw[line width=1.5pt, loosely dotted, gray](24)to(26);
					\draw[line width=1.5pt, loosely dotted, gray](25)to(27);
					\draw[line width=1.5pt, loosely dotted, gray](22)to(28);
					
					\draw[rounded corners=10pt,fill=gray, opacity=0.2] (-2.8,1.3) -- (-1.6,2.5) -- (0.1,2.5) -- (2.9,-0.15) -- (1.6,-1.2) -- (-2.9,-1.2) -- (-4.2,-.15) -- (-2.8,1.3);
					
					
					\node at (-3.5,-.4) [scale=1.25]{$v_1$};
					\node at (-2.5,-.4) [scale=1.25,fill=mygray]{\phantom{..}};
					\node at (-2.5,-.4) [scale=1.25]{$v_2$};
					\node at (-1.5,-.4) [scale=1.25,fill=mygray]{\phantom{..}};
					\node at (-1.5,-.4) [scale=1.25]{$v_3$};
					\node at (-.2,-.4) [scale=1.25,fill=mygray]{\phantom{..}};
					\node at (-.2,-.4) [scale=1.25]{$v_{2k-3}$};
					\node at (1.95,-.4) [scale=1.25, fill=mygray]{\phantom{-}};
					\node at (1.95,-.4) [scale=1.25]{$v_{2k-1}$};
					
					\draw[decorate, yshift=-4ex] (2.4,-.8) -- node[below=0.5ex] {~ } (-3.7,-.8);
					\node at (-.75,-2) [scale=1.5]{$2k-1$};
					\node at (-3.5,1.2) [scale=1.5]{$V_1$};
				\end{tikzpicture}
				\caption{We remove the grey vertices of the largest burning circle $V_1$ in a $2$-caterpillar.}
				\label{(3)}
			\end{center}
		\end{figure}
	\vspace{-.4cm}
		\noindent
		Analogously, for $\sum_{i=1}^{2k-2}\ell_\Sigma(v_i)= 1$ and $\ell_{\max}(v_{2k-1})\leq 1$ but  $\ell_{\max}(v_{2k-2})\geq 1$, we remove $v_1,\dots, v_{2k-2}$ with all adjacent legs. 
		Hence, it remains to consider the cases 
		\begin{enumerate}
			\item [a)] $\sum\limits_{i=1}^{2k-3}\ell_\Sigma(v_i)= 1$ with $\ell_{\max}(v_{2k-2})= 2$ and
			\item [b)] $\sum\limits_{i=1}^{2k-3}\ell_\Sigma(v_i)= 0$ with $\ell_{\max}(v_{2k-2})=2$ or $\ell_{\max}(v_{2k-1})\geq 1$.
		\end{enumerate} 
		\noindent
		If in case a) it additionally holds $$ \sum_{i=2k-1}^{(2k-1)+(2k-3)-4}\ell_\Sigma(v_i) \geq 1\quad\text{ or }\quad \ell_{\max}\left(v_{(2k-1)+(2k-3)-3}\right)\leq 1,$$ we arrange the two largest burning circles $V_1$ and $V_2$ with an overlap of two vertices as outlined in Figure \ref{(4)}. 
		We delete the vertices $v_1,\dots, v_{(2k-1)+(2k-3)-4}$ and, if $$\ell_{\max}\left(v_{(2k-1)+(2k-3)-3}\right)\leq 1,$$ we also remove $v_{(2k-1)+(2k-3)-3}$ with all its adjacent legs. Thus, at most $n-(2k-1)-(2k-3)$ vertices are left.
		
		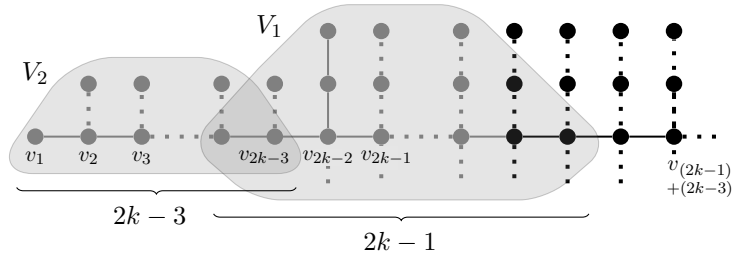
\begin{figure}[H]
			\begin{center}
				\begin{tikzpicture}[decoration=brace,node distance=5em, every node/.style={scale=0.6},scale=0.7]
				\node(1) at (-3.5,0)[circle,fill=gray]{};
				\node(2) at (1,0)[circle,fill=gray] {};
				\node(3) at (2,0)[circle,fill=gray] {};
				\node(4) at (4.5,0)[circle,fill=black] {};
				\node(5) at (5.5,0)[circle,fill=black] {};
				\node(6) at (6.5,0) [circle,fill= black] {};
				\node(7) at (2,1)[circle,fill=gray]{};
				\node(8) at (2,-1) {};
				\node(9) at (4.5,1)[circle,fill=black]{};
				\node(10) at (4.5,-1){};
				
				\node(11) at (-4.5,0) [circle,fill=gray] {};
				\node(12) at (-3.5,1)[circle,fill=gray]{};
				\node(14) at (1,1)[circle,fill=gray] {};
				\node(15) at (1,-1) {};
				\node(16) at (5.5,1)[circle,fill=black] {};
				\node(17) at (5.5,-1) {};
				\node(18) at (6.5,1)[circle,fill=black]{};
				\node(19) at (6.5,-1){};
				\node(20) at (7.5,0)[circle,fill=black] {};
				\node(21) at (0,0) [circle,fill=gray] {};
				\node(22) at (0,1)[circle,fill=gray]{};
				\node(24) at (-2.5,0) [circle,fill=gray] {};
				\node(25) at (-2.5,1)[circle,fill=gray]{};
				
				\node(29) at (1,2)[circle,fill=gray] {};
				\node(30) at (2,2)[circle,fill=gray] {};
				\node(31) at (-1,0)[circle,fill=gray] {};
				\node(32) at (-1,1)[circle,fill=gray] {};
				\node(33) at (3.5,0)[circle,fill=gray] {};
				\node(34) at (5.5,2)[circle,fill=black] {};
				\node(35) at (4.5,2)[circle,fill=black] {};
				\node(36) at (7.5,1)[circle,fill=black]{};
				\node(37) at (7.5,-1){};
				\node(38) at (8.5,0) {};
				\node(39) at (3.5,1)[circle,fill=gray]{};
				\node(40) at (3.5,-1) {};
				\node(41) at (3.5,2)[circle,fill=gray] {};
				\node(42) at (7.5,2)[circle,fill=black]{};
				\node(43) at (6.5,2)[circle,fill=black]{};
				
				\draw[line width=1.5pt, loosely dotted, gray](24)to(31);
				\draw[line width=.75pt, gray](2)to(3);
				\draw[line width=.75pt, gray](1)to(24);
				\draw[line width=1.5pt, loosely dotted, gray](3)to(33);
				\draw[line width=.75pt](4)to(5);
				\draw[line width=.75pt](6)to(20);
				\draw[line width=.75pt](5)to(6);
				\draw[line width=1.5pt, loosely dotted, gray](3)to(7);
				\draw[line width=1.5pt, loosely dotted, gray](3)to(8);
				\draw[line width=.75pt, gray](14)to(29);
				\draw[line width=1.5pt, loosely dotted, gray](7)to(30);
				\draw[line width=1.5pt, loosely dotted](4)to(9);
				\draw[line width=1.5pt, loosely dotted](4)to(10);
				\draw[line width=.75pt, gray](11)to(1);
				\draw[line width=1.5pt, loosely dotted, gray](12)to(1);
				\draw[line width=.75pt, gray](14)to(2);
				\draw[line width=1.5pt, loosely dotted, gray](15)to(2);
				\draw[line width=1.5pt, loosely dotted](16)to(5);
				\draw[line width=1.5pt, loosely dotted](17)to(5);
				\draw[line width=1.5pt, loosely dotted](18)to(6);
				\draw[line width=1.5pt, loosely dotted](19)to(6);
				\draw[line width=1.5pt, loosely dotted](20)to(38);
				\draw[line width=1.5pt, loosely dotted](35)to(9);
				\draw[line width=1.5pt, loosely dotted](34)to(16);
				\draw[line width=1.5pt, loosely dotted](36)to(20);
				\draw[line width=1.5pt, loosely dotted](37)to(20);
				\draw[line width=1.5pt, loosely dotted, gray](21)to(22);
				\draw[line width=1.5pt, loosely dotted, gray](24)to(25);
				\draw[line width=1.5pt, loosely dotted, gray](39)to(41);
				\draw[line width=1.5pt, loosely dotted, gray](33)to(39);
				\draw[line width=1.5pt, loosely dotted, gray](33)to(40);
				\draw[line width=1.5pt, loosely dotted, gray](31)to(32);
				\draw[line width=.75pt, gray](31)to(2);
				\draw[line width=.75pt, gray](33)to(4);
				\draw[line width=1.5pt, loosely dotted](20)to(42);
				\draw[line width=1.5pt, loosely dotted](18)to(43);
				
				\draw[rounded corners=15pt,fill=gray, opacity=0.2] (-4.55,.4) -- (-3.8,1.5) -- (-.7,1.5) -- (.8,-.7)  -- (-5.3,-.7) -- (-4.55,.4);
				\draw[rounded corners=10pt,fill=gray, opacity=0.2] (-.3,1.3) -- (0.9,2.5) -- (3.6,2.5) -- (6.3,-0.15) -- (5.1,-1.2) -- (-.4,-1.2) -- (-1.6,-.15) -- (-.3,1.3);
				
				
				\node at (-4.5,-.4) [scale=1.25]{$v_1$};
				\node at (-3.5,-.4) [scale=1.25]{$v_2$};
				\node at (-2.5,-.4) [scale=1.25]{$v_3$};
				\node at (1,-.4) [scale=1.25,fill=mygray]{\phantom{..}};
				\node at (1,-.4) [scale=1.25]{$v_{2k-2}$};
				\node at (-.2,-.4) [scale=1.25]{$v_{2k-3}$};
				\node at (2.1,-.4) [scale=1.25, fill=mygray]{\phantom{-}};
				\node at (2.1,-.4) [scale=1.25]{$v_{2k-1}$};
				\node at (7.7,-.7) [scale=1.25, fill=white]{\phantom{--}};
				\node at (8.0,-.8) [scale=1.25]{$v_{\substack{(2k-1)\\\hspace{-.2cm}+(2k-3)}}$};
				
				\draw[decorate, yshift=-4ex] (0.4,-.3) -- node[below=0.5ex] {~ } (-4.85,-.3);
				\node at (-2.35,-1.5) [scale=1.5]{$2k-3$};
				\node at (-4.5,1.2) [scale=1.5]{$V_2$};
				
				\draw[decorate, yshift=-4ex] (5.9,-.8) -- node[below=0.5ex] {~ } (-1.15,-.8);
				\node at (2.35,-2) [scale=1.5]{$2k-1$};
				\node at (-.1,2.1) [scale=1.5]{$V_1$};
				\end{tikzpicture}
				\caption{We arrange the two largest burning circles $V_1$ and $V_2$ with an overlap of two vertices.}
				\label{(4)}
			\end{center}
		\end{figure}
	\vspace{-.4cm}
		\noindent
		If, however, in case a) we have additionally $$\sum\limits_{i=2k-1}^{(2k-1)+(2k-3)-4}\ell_\Sigma(v_i) =0\quad\text{ and }\quad\ell_{\max}\left(v_{(2k-1)+(2k-3)-3}\right)=2,$$ we consider the three largest burning circles and position them as shown in Figure \ref{(5)}. 
		The removal of $v_1,\dots, v_{(2k-1)+(2k-3)-4}$ with all adjacent legs yields a graph with at most $n-(2k-1)-(2k-3)-(2k-5)-1$ vertices.
		
		\begin{figure}[H]
			\begin{center}
				\begin{tikzpicture}[decoration=brace,node distance=5em, every node/.style={scale=0.6},scale=0.7]
				\node(1) at (-3.5,0)[circle,fill=gray]{};
				\node(2) at (1,0)[circle,fill=gray] {};
				\node(3) at (2,0)[circle,fill=gray] {};
				\node(4) at (4.5+5.5,0)[circle,fill=black] {};
				\node(5) at (5.5+5.5,0)[circle,fill=black] {};
				\node(6) at (6.5+5.5,0) [circle,fill= black] {};
				\node(8) at (2,-1) {};
				\node(9) at (4.5+5.5,1)[circle,fill=black]{};
				\node(10) at (4.5+5.5,-1){};
				
				\node(11) at (-4.5,0) [circle,fill=gray] {};
				\node(12) at (-3.5,1)[circle,fill=gray]{};
				\node(14) at (1,1)[circle,fill=gray] {};
				\node(15) at (1,-1) {};
				\node(16) at (5.5+5.5,1)[circle,fill=black] {};
				\node(17) at (5.5+5.5,-1) {};
				\node(18) at (6.5+5.5,1)[circle,fill=black]{};
				\node(19) at (6.5+5.5,-1){};
				\node(20) at (7.5+5.5,0)[circle,fill=black] {};
				\node(21) at (0,0) [circle,fill=gray] {};
				\node(22) at (0,1)[circle,fill=gray]{};
				\node(24) at (-2.5,0) [circle,fill=gray] {};
				\node(25) at (-2.5,1)[circle,fill=gray]{};
				
				\node(29) at (1,2)[circle,fill=gray] {};
				\node(31) at (-1,0)[circle,fill=gray] {};
				\node(32) at (-1,1)[circle,fill=gray] {};
				\node(33) at (3.5+5.5,0)[circle,fill=gray] {};
				\node(34) at (5.5+5.5,2)[circle,fill=black] {};
				\node(35) at (4.5+5.5,2)[circle,fill=black] {};
				\node(36) at (7.5+5.5,1)[circle,fill=black]{};
				\node(37) at (7.5+5.5,-1){};
				\node(38) at (8.5+5.5,0) {};
				\node(39) at (3.5+5.5,1)[circle,fill=gray]{};
				\node(40) at (3.5+5.5,-1) {};
				\node(41) at (3.5+5.5,2)[circle,fill=gray] {};
				\node(42) at (7.5+5.5,2)[circle,fill=black]{};
				\node(43) at (6.5+5.5,2)[circle,fill=black]{};
				\node(44) at (3.5,0)[circle,fill=gray]{};
				\node(45) at (4.5,0)[circle,fill=gray]{};
				\node(46) at (5.5,0)[circle,fill=gray]{};
				\node(47) at (6.5,0)[circle,fill=gray]{};
				\node(48) at (7.5,0)[circle,fill=gray]{};
				\node(49) at (7.5,1)[circle,fill=gray]{};
				\node(50) at (7.5,2)[circle,fill=gray]{};

				\draw[line width=1.5pt, loosely dotted, gray](24)to(31);
				\draw[line width=.75pt, gray](2)to(3);
				\draw[line width=.75pt, gray](1)to(24);
				\draw[line width=1.5pt, loosely dotted, gray](3)to(44);
				\draw[line width=.75pt](4)to(5);
				\draw[line width=.75pt](6)to(20);
				\draw[line width=.75pt](5)to(6);
				\draw[line width=1.5pt, loosely dotted, gray](3)to(8);
				\draw[line width=.75pt, gray](14)to(29);
				\draw[line width=1.5pt, loosely dotted](4)to(9);
				\draw[line width=1.5pt, loosely dotted](4)to(10);
				\draw[line width=.75pt, gray](11)to(1);
				\draw[line width=1.5pt, loosely dotted, gray](12)to(1);
				\draw[line width=.75pt, gray](14)to(2);
				\draw[line width=1.5pt, loosely dotted, gray](15)to(2);
				\draw[line width=1.5pt, loosely dotted](16)to(5);
				\draw[line width=1.5pt, loosely dotted](17)to(5);
				\draw[line width=1.5pt, loosely dotted](18)to(6);
				\draw[line width=1.5pt, loosely dotted](19)to(6);
				\draw[line width=1.5pt, loosely dotted](20)to(38);
				\draw[line width=1.5pt, loosely dotted](35)to(9);
				\draw[line width=1.5pt, loosely dotted](20)to(42);
				\draw[line width=1.5pt, loosely dotted](18)to(43);
				\draw[line width=1.5pt, loosely dotted](34)to(16);
				\draw[line width=1.5pt, loosely dotted](36)to(20);
				\draw[line width=1.5pt, loosely dotted](37)to(20);
				\draw[line width=1.5pt, loosely dotted, gray](21)to(22);
				\draw[line width=1.5pt, loosely dotted, gray](24)to(25);
				\draw[line width=1.5pt, loosely dotted, gray](39)to(41);
				\draw[line width=1.5pt, loosely dotted, gray](33)to(39);
				\draw[line width=1.5pt, loosely dotted, gray](33)to(40);
				\draw[line width=1.5pt, loosely dotted, gray](31)to(32);
				\draw[line width=.75pt, gray](31)to(2);
				\draw[line width=.75pt, gray](33)to(4);
				\draw[line width=.75pt, gray](44)to(45);
				\draw[line width=.75pt, gray](46)to(45);
				\draw[line width=.75pt, gray](46)to(47);
				\draw[line width=.75pt, gray](47)to(48);
				\draw[line width=.75pt, gray](48)to(49);
				\draw[line width=.75pt, gray](49)to(50);
				\draw[line width=1.5pt, loosely dotted, gray](48)to(33);
				
				\draw[rounded corners=15pt,fill=gray, opacity=0.2] (-4.55,.4) -- (-3.8,1.5) -- (-.7,1.5) -- (.8,-.7)  -- (-5.3,-.7) -- (-4.55,.4);
				\draw[rounded corners=10pt,fill=gray, opacity=0.2] (-.3,1.3) -- (0.9,2.5) -- (2.6,2.5) -- (5.15,-0.15) -- (4.1,-1.2) -- (-.4,-1.2) -- (-1.6,-.15) -- (-.3,1.3);
				\draw[rounded corners=10pt,fill=gray, opacity=0.2] (-.3+6.5,1.3) -- (0.9+6.5,2.5) -- (2.6+6.5,2.5) -- (5.3+6.5,-0.15) -- (4.1+6.5,-1.2) -- (-.4+6.5,-1.2) -- (-1.65+6.5,-.15) -- (-.3+6.5,1.3);
				
				
				\node at (-4.5,-.4) [scale=1.25]{$v_1$};
				\node at (-3.5,-.4) [scale=1.25]{$v_2$};
				\node at (-2.5,-.4) [scale=1.25]{$v_3$};
				\node at (1,-.4) [scale=1.25,fill=mygray]{\phantom{..}};
				\node at (1,-.4) [scale=1.25]{$v_{2k-2}$};
				\node at (-.2,-.4) [scale=1.25]{$v_{2k-3}$};
				\node at (2.1,-.4) [scale=1.25, fill=mygray]{\phantom{-}};
				\node at (2.1,-.4) [scale=1.25]{$v_{2k-1}$};
				\node at (7.7+5.5,-.7) [scale=1.25, fill=white]{\phantom{--}};
				\node at (8.0+5.5,-.8) [scale=1.25]{$v_{\substack{(2k-1)\\\hspace{-.2cm}+(2k-3)}}$};
				\node at (7.7,-.4) [scale=1.25, fill=mygray]{\phantom{-}};
				\node at (7.7,-.7) [scale=1.25]{$v_{\substack{(2k-3)\\\hspace{-.19cm}+(2k-5)\\\hspace{-.96cm}+1}}$};
				
				\draw[decorate, yshift=-4ex] (0.4,-.3) -- node[below=0.5ex] {~ } (-4.85,-.3);
				\node at (-2.35,-1.5) [scale=1.5]{$2k-3$};
				\node at (-4.5,1.2) [scale=1.5]{$V_2$};
				
				\draw[decorate, yshift=-4ex] (4.9,-.8) -- node[below=0.5ex] {~ } (-1.15,-.8);
				\node at (1.85,-2) [scale=1.5]{$2k-5$};
				\node at (-.1,2.1) [scale=1.5]{$V_3$};
				
				\draw[decorate, yshift=-4ex] (4.9+6.5,-.8) -- node[below=0.5ex] {~ } (-1.15+6.5,-.8);
				\node at (1.85+6.5,-2) [scale=1.5]{$2k-1$};
				\node at (-.1+6.5,2.1) [scale=1.5]{$V_1$};
				\end{tikzpicture}
				\caption{We delete the grey vertices of the three largest burning circles $V_1$, $V_2$ and $V_3$.}
				\label{(5)}
			\end{center}
		\end{figure}
		\vspace{-.5cm}
		\noindent
		Lastly, in case b) we can assume without loss of generality that $\sum_{i=1}^{2k-3}\ell_\Sigma(v_i)= 0$ with $\ell_{\max}(v_{2k-2})=2$ or $\ell_{\max}(v_{2k-1})\geq 1$ holds for both ends of the spine (otherwise we can apply one of the cases above on the other end), i.e., additionally, we have $\sum_{i=1}^{2k-3}\ell_\Sigma(v_{l-i+1})= 0$. 
		Considering the three largest burning circles again, we place $V_3$ and $V_1$ at the beginning of the spine if $$\sum\limits_{i=2k-2}^{(2k-5)+(2k-1)-2}\ell_\Sigma(v_{i})\geq 2$$ 
		and at the end if 
		$$\sum\limits_{i=l-(2k-2)+1}^{l-((2k-5)+(2k-1)-2)+1}\ell_\Sigma(v_{i})\geq 2.$$ 
		As outlined in Figure \ref{(6)}, we put $V_2$ at the other side of the spine and remove the vertices $v_1,\dots, v_{(2k-5)+(2k-1)-2}$ and $v_l,\dots, v_{l-(2k-3)+1}$, respectively $v_{l},\dots, v_{l-((2k-5)+(2k-1)-2)+1}$ and $v_1,\dots, v_{2k-3}$.\\
		In the remaining case both sums equal one, $l_\Sigma(v_{2k-1})=l_\Sigma\left(v_{l-(2k-1)+1}\right)=1$ and $l_\Sigma(v_{i})=0$ for all other $(2k-5)+(2k-1)-2$ spine vertices at both ends. Thus, we incorporate $V_4$, placing it next to $V_2$ without overlap, and additionally remove $2k-7$ spine vertices, one of which has an adjacent leg.
		
		\begin{figure}[H]
			\begin{center}
				\begin{tikzpicture}[decoration=brace,node distance=5em, every node/.style={scale=0.6},scale=0.68]
				\node(1) at (-5.5,0)[circle,fill=gray]{};
				\node(2) at (1,0)[circle,fill=gray] {};
				\node(3) at (2,0)[circle,fill=gray] {};
				\node(4) at (4.5,0)[circle,fill=black] {};
				\node(5) at (5.5,0)[circle,fill=black] {};
				\node(6) at (6.5,0) [circle,fill= black] {};
				\node(7) at (2,1)[circle,fill=gray]{};
				\node(8) at (2,-1) {};
				\node(9) at (4.5,1)[circle,fill=black]{};
				\node(10) at (4.5,-1){};
				
				\node(14) at (1,1)[circle,fill=gray] {};
				\node(15) at (1,-1) {};
				\node(16) at (5.5,1)[circle,fill=black] {};
				\node(17) at (5.5,-1) {};
				\node(18) at (6.5,1)[circle,fill=black]{};
				\node(19) at (6.5,-1){};
				\node(20) at (7.5,0)[circle,fill=black] {};
				\node(24) at (-4.5,0) [circle,fill=gray] {};
				
				\node(29) at (1,2)[circle,fill=gray] {};
				\node(30) at (2,2)[circle,fill=gray] {};
				\node(31) at (-3,0)[circle,fill=gray] {};
				\node(33) at (3.5,0)[circle,fill=gray] {};
				\node(34) at (5.5,2)[circle,fill=black] {};
				\node(35) at (4.5,2)[circle,fill=black] {};
				\node(36) at (7.5,1)[circle,fill=black]{};
				\node(37) at (7.5,-1){};
				\node(38) at (9,0) {};
				\node(39) at (3.5,1)[circle,fill=gray]{};
				\node(40) at (3.5,-1) {};
				\node(41) at (3.5,2)[circle,fill=gray] {};
				\node(42) at (0,0)[circle,fill=gray] {};
				\node(43) at (-1,0)[circle,fill=gray] {};
				\node(44) at (-2,0)[circle,fill=gray] {};
				\node(45) at (7.5,2)[circle,fill=black]{};
				\node(46) at (6.5,2)[circle,fill=black]{};

				\draw[line width=1.5pt, loosely dotted, gray](24)to(31);
				\draw[line width=.75pt, gray](2)to(3);
				\draw[line width=.75pt, gray](1)to(24);
				\draw[line width=1.5pt, loosely dotted, gray](3)to(33);
				\draw[line width=.75pt](4)to(5);
				\draw[line width=.75pt](6)to(20);
				\draw[line width=.75pt](5)to(6);
				\draw[line width=1.5pt, loosely dotted,gray](3)to(7);
				\draw[line width=1.5pt, loosely dotted, gray](3)to(8);
				\draw[line width=1.5pt, loosely dotted, gray](14)to(29);
				\draw[line width=1.5pt, loosely dotted, gray](7)to(30);
				\draw[line width=1.5pt, loosely dotted](4)to(9);
				\draw[line width=1.5pt, loosely dotted](4)to(10);
				\draw[line width=.75pt, gray](11)to(1);
				\draw[line width=1.5pt, loosely dotted, gray](14)to(2);
				\draw[line width=1.5pt, loosely dotted, gray](15)to(2);
				\draw[line width=1.5pt, loosely dotted](16)to(5);
				\draw[line width=1.5pt, loosely dotted](17)to(5);
				\draw[line width=1.5pt, loosely dotted](18)to(6);
				\draw[line width=1.5pt, loosely dotted](19)to(6);
				\draw[line width=1.5pt, loosely dotted](35)to(9);
				\draw[line width=1.5pt, loosely dotted](34)to(16);
				\draw[line width=1.5pt, loosely dotted](36)to(20);
				\draw[line width=1.5pt, loosely dotted](37)to(20);
				\draw[line width=1.5pt, loosely dotted, gray](39)to(41);
				\draw[line width=1.5pt, loosely dotted, gray](33)to(39);
				\draw[line width=1.5pt, loosely dotted, gray](33)to(40);
				\draw[line width=.75pt, gray](31)to(2);
				\draw[line width=.75pt, gray](33)to(4);
				\draw[line width=1.5pt, loosely dotted](20)to(45);
				\draw[line width=1.5pt, loosely dotted](18)to(46);
				
				\draw[draw=gray, fill=gray, opacity=0.2] (-3.7,0) ellipse (2.25cm and 1cm);
				\draw[rounded corners=10pt,fill=gray, opacity=0.2] (-.3,1.3) -- (0.9,2.5) -- (3.6,2.5) -- (6.3,0) -- (5.1,-1.2) -- (-.4,-1.2) -- (-1.6,0) -- (-.3,1.3);
				
				
				\node at (-5.5,-.4) [scale=1.25]{$v_1$};
				\node at (-4.5,-.4) [scale=1.25]{$v_2$};
				\node at (-.1,-.4) [scale=1.25]{$v_{2k-3}$};
				\node at (-2.18,-.4) [scale=1.25]{$v_{2k-5}$};
				\node at (2.1,-.4) [scale=1.25, fill=mygray]{\phantom{-}};
				\node at (2.1,-.4) [scale=1.25]{$v_{2k-1}$};
				\node at (7.6,-.7) [scale=1.25, fill=white]{\phantom{--}};
				\node at (7.6,-.8) [scale=1.25]{$v_{\hspace{-.01cm}\substack{(2k-1)\\\hspace{-.14cm}+(2k-3)}}$};
				
				\draw[decorate, yshift=-4ex] (-1.6,-.5) -- node[below=0.5ex] {~ } (-5.85,-.5);
				\node at (-3.85,-1.7) [scale=1.5]{$2k-5$};
				\node at (-5.5,1.1) [scale=1.5]{$V_3$};
				
				\draw[decorate, yshift=-4ex] (5.9,-.8) -- node[below=0.5ex] {~ } (-1.15,-.8);
				\node at (2.35,-2) [scale=1.5]{$2k-1$};
				\node at (-.1,2.1) [scale=1.5]{$V_1$};
				

				\node(101) at (3.5+11.5,0)[circle,fill=gray]{};
				\node(102) at (-1+11.5,0)[circle,fill=black] {};
				\node(103) at (-2+11.5,0)[circle,fill=black] {};
				\node(107) at (-2+11.5,1)[circle,fill=black]{};
				\node(108) at (-2+11.5,-1) {};
				
				\node(114) at (-1+11.5,1)[circle,fill=black] {};
				\node(115) at (-1+11.5,-1) {};
				\node(121) at (0+11.5,0) [circle,fill=gray] {};
				\node(124) at (2.5+11.5,0) [circle,fill=gray] {};
				\node(129) at (-1+11.5,2)[circle,fill=black] {};
				\node(130) at (-2+11.5,2)[circle,fill=black] {};
				\node(131) at (1+11.5,0)[circle,fill=gray] {};
				\node(142) at (-3.5+11.5,0){};
				%
				
				\draw[line width=1.5pt, loosely dotted, gray](124)to(131);
				\draw[line width=.75pt, black](102)to(103);%
				\draw[line width=.75pt, gray](101)to(124);
				\draw[line width=1.5pt, loosely dotted, black](103)to(107);%
				\draw[line width=1.5pt, loosely dotted, black](103)to(108);
				\draw[line width=1.5pt, loosely dotted, black](114)to(129);
				\draw[line width=1.5pt, loosely dotted, black](107)to(130);
				\draw[line width=1.5pt, loosely dotted, black](114)to(102);
				\draw[line width=1.5pt, loosely dotted, black](115)to(102);
				\draw[line width=1.5pt, loosely dotted, black](103)to(20);
				\draw[line width=.75pt, gray](131)to(121);
				\draw[line width=.75pt, gray](102)to(121);
				
				\draw[draw=gray, fill=gray, opacity=0.2] (1.7+11.5,-.05) ellipse (2.25cm and 1cm);
				
				\node at (3.5+11.5,-.4) [scale=1.25]{$v_{l}$};
				\node at (2.5+11.5,-.4) [scale=1.25]{$v_{l-1}$};
				\node at (.45+11.5,-.6) [scale=1.25]{$v_{\substack{\hspace{0.03cm}l-(2k-3)\\+1}}$};
				\node at (-1.9+11.5,-.6) [scale=1.25, fill=white]{\phantom{H}};
				\node at (-1.9+11.5,-.7) [scale=1.25]{$v_{\substack{\hspace{0.03cm}l-(2k-1)\\+1}}$};
				
				\draw[decorate, yshift=-4ex] (3.85+11.5,-.5) -- node[below=0.5ex] {~ } (-0.5+11.5,-.5);
				\node at (1.85+11.5,-1.6) [scale=1.5]{$2k-3$};
				\node at (3.5+11.5,.9) [scale=1.5]{$V_2$};
				\end{tikzpicture}
				\caption{We delete the grey vertices of the three largest burning circles $V_1$, $V_2$ and $V_3$.}
				\label{(6)}
			\end{center}
		\end{figure}
		\vspace{-.5cm}
		\noindent
		This completes the proof of the Burning Number Conjecture for $2$-caterpillars.
		\vspace{-.55cm}
	\end{proof}

\end{theorem}
\vspace{.55cm}
\noindent
Finally, we prove the conjecture for all $p$-caterpillars with at least $2\left\lceil\sqrt{n}\right\rceil-1$ vertices of degree one, i.e., for all trees having at least $2\left\lceil\sqrt{n}\right\rceil-1$ leaves.

\begin{theorem}
	The burning number of a $p$-caterpillar $G$ with at least $2\left\lceil\sqrt{n}\right\rceil-1$ vertices of degree one satisfies $b(G)\leq \left\lceil\sqrt{n}\right\rceil$.
	\begin{proof}
		Assume $G=(V,E)$ to be a minimum counterexample regarding $p$ and among these minimal regarding the vertex number $n$. 
		Hence, $b(G)>\left\lceil\sqrt{n}\right\rceil=:k$ and $\vert L\vert\geq 2k-1$ with the notation $L:=\{v\in V~\vert~\deg(v)=1\}$. 
		Deleting all leaves, the remaining graph $G-L$ is  a $(p-1)$-caterpillar and thus, 
			\begin{align*}
				b(G-L)
				\leq&\left\lceil\sqrt{n-\vert L\vert}\right\rceil
				\leq\left\lceil\sqrt{n-2k+1}\right\rceil
				\leq\left\lceil\sqrt{n-2\sqrt{n}+1}\right\rceil\\
				=&\left\lceil\sqrt{(\sqrt{n}-1)^2}\right\rceil
				=\left\lceil\sqrt{n}\right\rceil-1.
			\end{align*}
		However, if $G-L$ burns after $\left\lceil\sqrt{n}\right\rceil-1$ steps, using the same burning strategy $G$ can be burnt in $\left\lceil\sqrt{n}\right\rceil$ steps. This contradicts the assumption and thus, no counterexample exists.
		\vspace{-.55cm}
	\end{proof}
\end{theorem}


\section{Concluding Remarks}\label{sec:Conclusion}


\noindent
By the results of this paper, it remains to prove the conjecture for $p$-caterpillars, $p\geq3$, with less than $2\left\lceil\sqrt{n}\right\rceil-1$ leaves to complete the proof of the conjectured bound for all connected graphs. 
Minimum counterexamples for these remaining graph classes can be characterised in great detail. 
We plan to investigate these characterisations to prove the conjecture in future work.

\bibliographystyle{abbrv} 
\bibliography{Literatur_BurningNumber}
\end{document}